\documentclass[12pt,twoside]{article}
\usepackage{amssymb,mathrsfs,amsmath,ifpdf,graphicx}

\long\def\@makefntext#1{\noindent #1}
\newskip\tabcentering \tabcentering=1000pt plus 1000pt minus 1000pt
\def\REF#1{\par\hangindent\parindent\indent\llap{#1\enspace}\ignorespaces} 
\def\MCH#1#2{\setbox0=\hbox{\raise#1\hbox{#2}}\smash{\box0}}

\def\@evenfoot{}\def\@oddfoot{}

\def\@evenhead{\hbox to\textwidth{\footnotesize\hfill\rm\thepage
}} 

\def\@oddhead{\hbox to \textwidth{\footnotesize \hfill\rm\thepage}}

\def\sec#1{\vskip 3mm\leftline{\bf #1}\vskip 1mm}

\def\bc{\begin{center}}
\def\ec{\end{center}}
\def\no{\noindent}
\def\hang{\hangindent\parindent}
\def\textindent#1{\indent\llap{\qquad #1\ \ \enspace}\ignorespaces}
\def\ref{\par\hang\textindent}

\begin{document}
\abovedisplayskip=6pt plus 1pt minus 1pt \belowdisplayskip=6pt
plus 1pt minus 1pt
\thispagestyle{empty} \vspace*{-1.0truecm} \noindent

\vskip 8mm

\bc{\large\bf The thickness of cartesian product $K_n \Box P_m$$^\ast$

\footnotetext{\footnotesize $^\ast$Supported by the National Natural Science Foundation of China under Grant No. 11126167}}\ec

\vskip 3mm \bc{\bf Yan Yang}\ec
\small \bc{Department of Mathematics}\vspace{-4mm}\ec
\small \bc{Tianjin University, Tianjin {\rm 300072,} P.R.China}\vspace{-4mm}\ec
\small \bc{E-mail: yanyang@tju.edu.cn}\ec

\vspace{3mm}
\noindent{\small {\small\bf Abstract} \ \ The thickness $\theta(G)$ of a graph $G$ is the minimum number of planar spanning subgraphs  into
which the graph $G$ can be decomposed. It is a topological invariant of a graph, which was defined by W.T. Tutte in $1963$ and
also has important applications to VLSI design. But comparing with other topological invariants, e.g. genus and crossing number, results about thickness of graphs are few.  The only types of graphs whose thicknesses have been obtained are complete graphs, complete bipartite graphs and hypercubes. In this paper, by operations on graphs, the thickness of the cartesian product $K_n \Box P_m$, $n,m \geq 2$ are obtained.

\vspace{2mm}
\no{\small\bf Keywords} \ \ thickness; cartesian product; complete graph; path.

\vspace{2mm}
\no{\small\bf Mathematics Subject Classifications} \ \ 05C10.

\vspace{5mm}
\sec{\large 1\quad Introduction}\vspace{1mm}

\no Let $G$ be a graph with  vertex set $V(G)$ and edge set $E(G)$, and graphs in this paper are all simple graphs. The  {\it cartesian product} of graphs $G$ and $H$ is the graph $G\Box H$ whose vertex set
$V(G\Box H)=V(G)\times V(H)$ and whose edge set $E(G\Box H)=\{(u_1,v_1)(u_2,v_2)~|~ u_1u_2\in E(G) ~\mbox{and}~ v_1=v_2, ~~\mbox {or}~~ v_1v_2\in E(H) ~\mbox{and}~ u_1=u_2\}$.

The {\it union }of graphs $G_1$ and $G_2$ is the graph $G_1\cup G_2$
with vertex set $V (G_1) \cup V (G_2)$ and edge set $E(G_1) \cup E(G_2)$.
The {\it intersection} $G_1 \cap G_2$ of $G_1$ and $G_2$ is defined analogously.
Let $G_1$ and $G_2$ be subgraphs of a graph $G$. If $G=G_1\cup G_2$ and $G_1 \cap G_2=\{v\}$ (a vertex of $G$),
then we say that $G$ is the {\it vertex-amalgamation} of $G_1$ and $G_2$ at vertex $v$, denoted  $G=G_1\vee_{\{v\}}G_2$.

A graph is said to be {\it planar} if it can be  drawn on the plane so that no two edges cross (i.e., its edges meet only at their common ends). A planar graph drawn in this way is called a {\it plane graph}. Suppose $G_1, G_2, \ldots, G_k$ are spanning subgraphs of $G$, if $$E(G_1)\cup E(G_2)\cup\cdots \cup E(G_k)=E(G) ~\mbox{and} ~ E(G_i)\cap E(G_j)=\emptyset,  (i\neq j,~i,j=1,2,\ldots,k),$$ then $\{G_1, G_2, \ldots, G_k\}$ is a decomposition of $G$. Furthermore, if $G_1, G_2, \ldots, G_k$ are all planar graph, then $\{G_1, G_2, \ldots, G_k\}$ is a planar decomposition of $G$. The minimum number of planar spanning subgraphs into
which a graph $G$ can be decomposed is called the {\it thickness} of $G$, denoted by $\theta(G)$.

The thickness of a graph first defined by W.T. Tutte[6] in $1963$, as a  topological invariant of a graph, it is an important research object in topological graph theory, and it also has important applications to VLSI design[1]. But the thickness problems are very difficult, the results about them are few. The only types of graphs whose thicknesses have been obtained are complete graphs[2], complete bipartite graphs[3] and hypercubes[4]. The reader is
referred to [5] for more background and results about the thickness problems.

In this paper,  the thickness of the cartesian product $K_n \Box P_m$, $n,m \geq 2$ is obtained, in which $K_n$ is the complete graph with $n$ vertices and $P_m$ is the path with $m$ vertices.

\vspace{5mm}

 \sec{\large 2\quad The lower bound for the thickness of  $K_n \Box P_m$, $n,m \geq 2$}\vspace{1mm}

\no
 From Euler's polyhedron formula, a planar graph $G$ has at most $3|V|-6$ edges, furthermore, one can get the lower bound for the thickness of a graph as follows.
\vskip2mm
 \no{\bf Lemma 2.1}[3]\quad{\it  If $G=(V,E)$ is a graph with $|V|=n$ and $|E|=m$, then $\theta(G)\geq\big \lceil \frac{m}{3n-6} \big \rceil$.}
\vskip2mm
 As for the graph $K_n \Box P_m$, the numbers of vertices and edges are $|V|=mn$ and $|E|=\frac{mn(n-1)}{2}+n(m-1)$ respectively. By using Lemma 2.1,
 \begin{eqnarray*}\theta(K_n \Box P_m)
&\geq &\Big \lceil\frac{\frac{mn(n-1)}{2}+n(m-1)}{3mn-6} \Big\rceil\\
&=&\Big \lceil\frac{n+1}{6}+\frac{2}{6(mn-2)} \Big \rceil\\
&=&\left\{\begin{array}{cc}
       p+1& ~~~\mbox {when}~ n=6p,6p+1,6p+2,6p+3,6p+4;\\
        p+2& ~~~\mbox {when}~ n=6p+5.
              \end{array}\right.~~~~~~~~~~~~(1)\end{eqnarray*}
\vskip2mm
\no{\bf Lemma 2.2}[2]\quad{\it The thickness of the complete graph $K_n$ is $\theta(K_n)=\big \lfloor\frac{n+7}{6}\big \rfloor$, except that $\theta(K_9)=\theta(K_{10})=3.$}

\vspace{5mm}

 \sec{\large 3\quad The thickness of  $K_n \Box P_2$, $n\geq 2$}\vspace{1mm}

\no Let $K^1_n$ be the complete graph with $n$ vertices labeled by $v_1,v_2,\ldots, v_n$ respectively,
$K^2_n$ is a copy of $K^1_n$ and it's vertices labeled by $u_1,u_2,\ldots, u_n$ respectively, by joining the vertices $v_i$ and $u_i$ by an edge $v_iu_i$, $1\leq i\leq n$, one can get the graph $K_n \Box P_2$. Figure 1 illustrates $K_5 \Box P_2$.

By inserting a vertex $w_i$ at edge $v_iu_i$, for $1\leq i\leq n$, and merging these $n$ $2$-valent vertices $w_1, w_2,\ldots,w_n$ into
one vertex $w$, one can get a new graph. This graph also can be seen as  the vertex-amalgamation of $K_{n+1}$ and $K_{n+1}$ at $w$, denote by $K_{n+1}\vee_{\{w\}}K_{n+1}$. Figure 2 shows the graph $K_6 \vee_{\{w\}}K_6$.

\vskip2mm
\input{fig1.TpX}
 \begin{center}{~~~~{\bf Figure 1} \ \ The graph $K_5 \Box P_2$}\end{center}

\vskip2mm
\input{fig2.TpX}
\begin{center}{~~~~{\bf Figure 2} \ \ The graph $K_6 \vee_{\{w\}}K_6$}\end{center}

\vskip2mm
 \no{\bf Lemma 3.1}[7]\quad{\it If $G$ is the vertex-amalgamation of $G_1$ and $G_2$, $\theta(G_1)=n_1$ and $\theta(G_2)=n_2$, then $\theta(G)= \max\{n_1, n_2\}.$}
\vskip2mm
From Lemma 3.1, the thickness of $K_{n+1}\vee_{\{w\}}K_{n+1}$ is the same as the thickness of $K_{n+1}$. Let $\theta(K_{n+1})=t$ and $\{G_1,G_2,\ldots,G_t\}$ is a planar decomposition of $K_{n+1}$, then one can get a planar decomposition of $K_{n+1}\vee_{\{w\}}K_{n+1}$ as follows,
$$\{G_1\vee_{\{w\}}G_1, G_2\vee_{\{w\}}G_2,\ldots, G_t\vee_{\{w\}}G_t\},$$ in which $G_i\vee_{\{w\}}G_i$, $1\leq i\leq t$ are plane graphs. In this way, a planar decomposition of $K_6 \vee_{\{w\}}K_6$ is shown in Figure 3.

\vskip2mm
\input{fig3.TpX}
\begin{center}{~~~~{\bf Figure 3} \ \ A planar decomposition of $K_6 \vee_{\{w\}}K_6$}\end{center}
\vskip2mm

 From the construction of $G_i\vee_{\{w\}}G_i$, if the edge $v_qw\in G_i\vee_{\{w\}}G_i$, then $u_qw\in G_i\vee_{\{w\}}G_i$, $1\leq q\leq n$.
 For each graph $G_i\vee_{\{w\}}G_i$, $1\leq i\leq t$, if $v_qw, ~ u_qw\in G_i\vee_{\{w\}}G_i$, then we replace them by a new edge $v_qu_q$, for $1\leq q\leq n$, and delete the vertex $w$. In this way, we can get a new planar decomposition, which is exactly a planar decomposition of $K_n \Box P_2$. Figure 4 illustrates a planar decomposition of  $K_5 \Box P_2$ by using this way.

\vskip2mm
\input{fig4.TpX}
\begin{center}{~~~~{\bf Figure 4} \ \ A planar decomposition of $K_5 \Box P_2$}\end{center}
\vskip2mm

From the argument and construction the above, one can get a planar decomposition of $K_n \Box P_2$ from that of $K_{n+1}\vee_{\{w\}}K_{n+1}$, so we have that $$~~~~~~~~~~~~~~~~~~~~~~~~~~~~\theta(K_n \Box P_2)\leq \theta(K_{n+1}\vee_{\{w\}}K_{n+1})=\theta(K_{n+1}).~~~~~~~~~~~~~~~~~~~~~~~(2)$$

\vskip2mm
 \no{\bf Theorem 3.1}\quad{\it  The thickness of the cartesian product $K_n \Box P_2$$(n\geq 2)$ is $\theta(K_n \Box P_2)=\big \lfloor\frac{n+8}{6}\big \rfloor$,
  except that $\theta(K_8 \Box P_2)=\theta(K_9 \Box P_2)=3$ and possibly when $n=6p+4$} ($n$ {\it is a nonnegative integer}).
 \vskip1mm
{\it Proof}\quad When $n\neq8,9$, from $(1),(2)$ and Lemma 2.2,  we can get that $\theta(K_n \Box P_2)=\theta(K_{n+1})$,
  except possibly when $n=6p+4$ ($p$ is a nonnegative integer).

When $n=8$, we suppose that $\theta(K_8 \Box P_2)=2$, and $\{G_1, G_2\}$ is a planar decomposition of $K_8 \Box P_2$ in which $G_1$ and $G_2$ are plane graphs. We insert $2$-valent vertex $w_i$ at edge $v_iu_i$, for each $1\leq i\leq n$, in $G_1$ and $G_2$, then merge these $2$-valent vertices into
vertices $w$ and $\tilde{w}$ in $G_1$ and $G_2$ respectively. In this way, we get two new plane graphs, denote by $\tilde{G_1}$ and $\tilde{G_2}$, and $\{\tilde{G_1}, \tilde{G_2}\}$ is a planar decomposition of $K_{9}\vee_{\{w\}}K_{9}$. From $(2)$, $\theta(K_{9})=\theta(K_{9}\vee_{\{w\}}K_{9})\leq 2$, it is a contradiction to the fact $\theta(K_{9})=3$. Hence, $\theta(K_8 \Box P_2)=\theta(K_{9})=3$. With a similar argument, we can get $\theta(K_9 \Box P_2)=\theta(K_{10})=3$.

Summarizing the above,  the theorem is
obtained.\hspace{\fill}$\Box$

\vspace{5mm}

 \sec{\large 4\quad The thickness of  $K_n \Box P_m$, $n\geq 2, m\geq 3$}\vspace{1mm}

\no We use the similar method with that in section 3.  Firstly, we insert a $2$-valent vertex into each "path edge" (the edges come from $P_m$).
Secondly, we merge these $(m-1)n$ $2$-valent vertices into $m-1$ vertices, each of which joint two adjacent $K_n$,  then we get a new graph $\tilde{G}$. The graph $\tilde{G}$ can be seen as a vertex-amalgamation of $m$ graphs, in which the first and the $m$th graphs are $K_{n+1}$, the others are $K_{n+2}-e$. From Lemma 3.1, one can get that
$$~~~~~~~~~~~~~~~~~~~~~~~~~~~~~~~~\theta(\tilde{G})=\theta(K_{n+2}-e)\leq \theta(K_{n+2}).~~~~~~~~~~~~~~~~~~~~~~~~~~~~~~~~~~~~~~~(3)$$

With a similar argument to that in section 3, we can get a planar decomposition of $K_n \Box P_m$ from a planar decomposition of $\tilde{G}$,  hence
$$~~~~~~~~~~~~~~~~~~~~~~~~~~~~~~\theta(K_n \Box P_m)\leq \theta(K_{n+2}-e)\leq \theta(K_{n+2}).~~~~~~~~~~~~~~~~~~~~~~~~~~~~~~~~~~(4)$$

\vskip1mm
 \no{\bf Theorem 4.1}\quad{\it  The thickness of the cartesian product $K_n \Box P_m$$(n\geq 2, m\geq 3)$ is \\ $\theta(K_n \Box P_m)=\big \lfloor\frac{n+9}{6}\big \rfloor$,
  except possibly when $n=6p+3, 6p+4$ and $n=8$} ($p$ {\it is a nonnegative integer}).
 \vskip1mm
{\it Proof}\quad When $n\neq 7,8$, from $(1),(4)$ and Lemma 2.2, we can get that $\theta(K_n \Box P_m)=\theta(K_{n+2})$,
  except possibly when $n=6p+3, 6p+4$ ($p$ is a nonnegative integer).

When $n=7$, from $(1)$ and $(4)$, we have $2\leq \theta(K_7 \Box P_m)\leq \theta(K_{n+2}-e)$.  We
can give a  planar decomposition of $K_{9}-e$ as shown in  Figure 5, and $K_{9}-e$ is a non-planar graph, hence $\theta(K_{9}-e)=2$, $\theta(K_7 \Box P_m)=2.$

Summarizing the above,  the theorem follows.\hspace{\fill}$\Box$

\input{fig5.TpX}
 \vspace{-0.3cm}\begin{center}{~~~~{\bf Figure 5} \ \ A planar decomposition of $K_9-e$}\end{center}

 \sec{\large 5\quad Concluded Remarks }\vspace{1mm}

\no Determining the thickness for a arbitrary graph is a difficult problem. The known results and methods are few. In [7], the authors first considered getting the thickness of graphs by operations on graphs. In this paper, the author use operations on graphs and some conclusions in [7], obtain the thickness of the cartesian product $K_n \Box P_m$.

For the exceptional cases in theorem 3.1,  when $n=6p+4$, $\theta(K_n \Box P_2)$ is either $\lfloor\frac{n+8}{6}\rfloor$ or $\lfloor\frac{n+8}{6}\rfloor-1$. For the exceptional cases in theorem 4.1,  when $n=6p+3, 6p+4$, $\theta(K_n \Box P_m)$ is either $\lfloor\frac{n+9}{6}\rfloor$ or $\lfloor\frac{n+9}{6}\rfloor-1$; when $n=8$, $\theta(K_8 \Box P_m)$ is either $2$ or $3$.
 What are the exactly numbers for these exceptional cases are still open.

\no \vskip0.2in \no {\bf References}

\footnotesize \

\REF{[1]} A. Aggarwal, M.Klawe, and P.Shor, Multilayer grid embeddings for VLSI, Algorithmica, 6(1991) 129-151го

\REF{[2]} L.W. Beineke and F. Harary, The thickness of the complete graph, Canad. J. Math., 17 (1965) 850-859го

\REF{[3]} L.W. Beineke, F. Harary, and J.W. Moon, On the thickness of the complete bipartite graph, Proc. Cambridge. Philo. Soc., 60 (1964) 1-5го

\REF{[4]} M. Kleinert,  Die Dicke des n-dimensionale W$\ddot{\mbox{u}}$rfel-Graphen, J. Comb. Theory, 3 (1967) 10-15.

\REF{[5]} E. M$\ddot{\mbox{a}}$kinen and T. Poranen. An annotated bibliography on the thickness, outerthickness, and arboricity of a graph. Technical report, University of Tampere, 2009.
\\http://www.sis.uta.fi/~tp54752/pub/thickness-bibliography/.

\REF{[6]} W.T. Tutte, The thickness of a graph, Indag. Math., 25 (1963) 567-577го

\REF{[7]} Y. Yang, X.H. Kong. The thickness of amalgamations of graphs. http://arxiv.org/abs/1201.6483.

\end{document}